\documentclass[12pt]{article}
\makeatletter
        \def\@thefnmark{\null}
        \def\footnotetexta{\@footnotetext}

\title{The diagonalizable nonnegative inverse eigenvalue problem}

\author{Anthony G. Cronin\footnote{School of Mathematics and Statistics, University College Dublin, Ireland, anthony.cronin@ucd.ie} and  Thomas J. Laffey\footnote{School of Mathematics and Statistics, University College Dublin, Ireland,  thomas.laffey@ucd.ie } }

\date{}
\usepackage{amsmath,amsfonts,mathrsfs,amssymb,amsthm}
\newtheorem{Theo}{Theorem}
\newtheorem{Lem}[Theo]{Lemma}

\newtheorem{Pro}[Theo]{Proposition}
\newtheorem{Rem*}[Theo]{Remark}

\newcommand{\npmatrix}[1]{\begin{bmatrix} #1 \end{bmatrix}}

\usepackage{enumerate}

\usepackage{color,epsfig,graphicx}
\usepackage{setspace}
\def\bbbc{\mathpalette{}{{\setbox0=\hbox{$\rm C$}\hbox{\hbox
to0pt{\kern0.4\wd0\vrule height0.9\ht0\hss}\box0}}}}

\begin{document}
\maketitle

\begin{abstract}
\noindent
In this paper we prove that the SNIEP $\neq$ DNIEP, i.e. the symmetric and diagonalizable nonnegative inverse eigenvalue problems are different. We also show that the minimum $t>0$ for which $(3+t,3-t,-2,-2,-2)$ is realizable by a diagonalizable matrix is $t=1$, and we distinguish diagonalizably realziable lists from general realizable lists using the Jordan Normal Form.
\end{abstract}

\section{Introduction}
Classifying spectra of nonnegative matrices is know as the nonnegative inverse eigenvalue problem or NIEP. The real nonnegative inverse eigenvalue problem or RNIEP is to determine necessary and sufficient conditions on the list of $n$ real numbers $(\lambda_1, \lambda_2, \ldots,\lambda_n)$ so that $\sigma$ be the spectrum of an entry-wise nonnegative matrix. Such a matrix is then said to \emph{realize} $\sigma$. If we further require that the realizing matrix be symmetric we call the problem the symmetric nonnegaitive inverse eigenvalue problem or SNIEP. Both problems are unsolved for $n\geq 5$ though the trace zero case for $5\times 5$ matrices was solved by Spector \cite{Spec}. The two problems are the same for $n\leq 4$ but are different for $n\geq5$ as was first shown by Johnson, Laffey and Loewy in \cite{JLL2}.
The first significant breakthrough in this problem was Suleimanova's result \cite{Sul} on lists of real numbers with just one positive number, which says that the real list $\lambda_1>0\geq \lambda_2\geq \cdots \geq \lambda_n$ is realizable if and only if $\lambda_1+ \lambda_2+ \cdots + \lambda_n \geq 0$. The question of realizing real lists of five or more numbers containing just two positive numbers is still unsolved in general. 
\subsection{Necessary conditions}
The Perron-Frobenius theory (\cite{Per},\cite{Fro}) says (among other things) that the spectral radius of an irreducible nonnegative matrix must be contained in the spectrum of that matrix i.e.
$$max\{|\lambda|_j:\lambda_j \in \sigma\}\in \sigma.$$
In this context the spectral radius $\rho$ is known as the \emph{Perron root} and for irreducible nonnegative matrices $\rho>0$ and it occurs just once as an eigenvalue.
We define the power sums $s_k$ as follows
$$s_k:=\lambda_1^k+\lambda_2^k+\cdots+\lambda_n^k, \mbox{ \ for \ } k=1,2,\ldots.$$
Notice that if $\sigma=(\lambda_1, \lambda_2, \ldots,\lambda_n)$ is the spectrum of a nonnegative matrix $A$ then the power sum $s_k$ is also the trace of the $k^{th}$ power of a realizing matrix $A$ for $\sigma$. Independently Loewy and London \cite{Loe/Lon} and Johnson \cite{JLL2} derived an infinite set of inequalities which the spectrum of a nonnegative matrix must satisfy, namely that 
$$n^{m-1}s_{km}\geq s_k^m \ \mbox{for} \ k,m=1,2, \ldots$$ known as the JLL conditions.
Necessary conditions for both the RNIEP and the SNIEP thus include:

\begin{eqnarray}\label{Perron}
\mbox{max}\{|\lambda|_j:\lambda_j \in \sigma\}\in \sigma
\end{eqnarray}
\begin{eqnarray}\label{Trace}
s_k\geq0 \mbox{ \ for \ } k=1,2, \ldots 
\end{eqnarray}
\begin{eqnarray}\label{JLL}
n^{m-1}s_{km}\geq s_k^m \mbox{ \ for \ } k,m,n=1,2, \ldots.
\end{eqnarray}
A new necessary condition for the SNIEP when the trace is at least half the spectral radius is given by Loewy and Spector in \cite{SpecL}.
A necessary condition for general $n$ involving only the first three power sums $s_k$ is given by Cronin and Laffey in \cite{CroLaf}.
\section{A classic example $\sigma=\{3,3,-2,-2,-2\}$}
The list $\sigma=\{3,3,-2,-2,-2\}$, in the guise $\tau=(1,1,-\frac{2}{3},-\frac{2}{3},-\frac{2}{3})$ was first studied by Salzmann in 1971 \cite{Sal} and Friedland in 1977 \cite{Fri}.
As can be checked the list $\tau$ satisfies the necessary conditions \eqref{Perron}, \eqref{Trace} and \eqref{JLL} for all positive integers $k,m$ and $n$.\\
However the list $\tau$ is not realizable. For if it were, a realizing matrix would have to be reducible as the Perron root occurs twice. As can be easily seen, any such reducible matrix having this spectrum would require the list to be partionable into two lists which are separately realizable, and this leads to a negative trace \emph{realization} which is forbidden by \eqref{Trace}.\\
Laffey and Meehan in \cite{LafMee1} showed that in order for a list of five numbers which sum to zero to be realizable, we must have a refined JLL inequality satisfied, in particular $4s_4-s_2^2 \geq 0$. For $\sigma_t$ (which is $3\tau$) we have $4s_4-s_2^2=840-900<0$ and so again we see that $\sigma$ cannot be realizable. \\

\subsection{$\sigma_t=\{3+t,3,-2,-2,-2\}$}
The related problem of trying to realize $\sigma_t=\{3+t,3,-2,-2,-2\}$ for the minimum value of $t>0$ is unsolved.\\

In her thesis, Meehan \cite{Mee} showed that $\sigma_t$ is realizable for $t\geq 0.519310982048\cdots$ and a realizing matrix of the form
$$A=\npmatrix{ t & 1 & 0 & 0 & 0  \cr p & 0 & 1 & 0 & 0 \cr
0 & q & 0 & 1 & 0 \cr 0 & 0 & 0 & 0 & 1 \cr 0 &
0 & w & h & 0 \cr}$$ is presented.
She also shows that for a $5 \times 5$ extreme (or Perron extreme) matrix $\cite{LafExt}$, \\ 
$4s_4-s_2^2+s_1^2s_2-\frac{s_1^4}{2}\geq 0$ which for $\sigma_t$, means $t\geq 0.39671\cdots$.\\
Hence the range for the minimum value of $t$ for which $\sigma_t$ is realizable is
$$0.39671\cdots \leq t \leq 0.51931\cdots$$
This classic example highlights the difficulty in moving from the NIEP for $n=5$ and trace zero to the positive trace case in the NIEP, even when the numbers of the list are all real. 
\subsubsection{$\sigma_t=\{3+t,3,-2,-2,-2\}$ in the symmetric case}
We note that for $t=1$, $\sigma_t=(3+t,3,-2,-2,-2)$, is \emph{symmetrically} realizable by the matrix 
$$A=\npmatrix{ 0 & 2 & 2 & 0 & 0  \cr 2 & 0 & 2 & 0 & 0 \cr
2 & 2 & 0 & 0 & 0 \cr 0 & 0 & 0 & 1 & \sqrt{6} \cr 0 &
0 & 0 & \sqrt{6} & 0 \cr}$$
and this is the best possible $t$ in the symmetric case.\\
Thus the RNIEP and the SNIEP are different for $n=5$, and this is the first case in which they differ $\cite{JLL2}$. The interested reader should consult Loewy and London \cite{Loe/Lon}, for the proofs that the RNIEP=SNIEP for $n \leq 4$.


\subsection{$\widehat{\sigma}=\{3+t,3-t,-2,-2,-2\}$}
A further related problem of finding the minimum value of $t>0$ for which the list \\
$\widehat{\sigma_t}=(3+t,3-t,-2,-2,-2)$ is realizable, is solved. Note that the sum of the elements in $\widehat{\sigma_t}$ is zero and so any realizing matrix for $\widehat{\sigma_t}$ must have trace zero. Also note, as is the case for $\sigma=(3+t,3-t,-2,-2,-2)$, that any realizing matrix for $\sigma$ must be irreducible.\\
The minimum value of $t$ for which $\widehat{\sigma_t}=(3+t,3-t,-2,-2,-2)$ is realziable in the \emph{symmetric} case must be at least one. This can be seen by considering the sign pattern of the eigenvector $z$ associated with the second largest eigenvalue $3-t$ of a realizing matrix $A=A^t$. Now $z$ cannot have all its entries positive (or negative) as only the Perron eigenvector (and its scalar multiples) has this property. Thus $z$ (or$-z$) has either one or two positive entries (and hence either four or three non-positive entries). Without loss of generality we can permute the entries of $z$ so that its positive entries, denoted with a + below,  occur in the first position(s) of the vector i.e.
$$z=\npmatrix{ +  \cr  - \cr - \cr - \cr -  } \mbox{ \ or \ } \npmatrix{ + \cr + \cr - \cr - \cr - }$$
where the minus sign means the corresponding entry is less than or equal to zero.\\
Now if $z$ has just one positive entry then $Az=(3-t)z$ implies
$$\npmatrix{ 0  & A_{12}\cr   
 A_{21} &  A_{22}   }\left[\begin{array}{c} z_1\\   z_2 \end{array}\right] 
 =(3-t)\left[\begin{array}{c} z_1\\   z_2 \end{array}\right]  $$
where $A_{12}$ is $(1 \times 4)$, $A_{21}$ is $(4 \times1)$, and $A_{22}$ is $(4 \times 4)$, and $z=[z_1 \ z_2]^t$ where $z_1>0$ and $z_2\leq 0$, for $z_1 \in \mathbb{R}, z_2 \in \mathbb{R}^4$. But this implies $A_{12}z_2^t=(3-t)z_1$ which is false as the left hand side of this equation is non-positive and the right hand side is positive. Hence this sign pattern (one positive component followed by four non-positive) for $z$ is not permitted. So $z$ must have 2 positive entries and 3 non-positive entries. \\
For simplicity we write $z=[u \ -v]^t$ where $u=[u_1 \ u_2]^t$ and $v=[v_1 \ v_2 \ v_3]^t$ where $u_i,v_j \geq 0$ for each $i$ and $j$. From the eigenvalue/eigenvector equation $Az=(3-t)z$ we have that 
$$\npmatrix{ A_{11}  & A_{12}\cr   
 A_{21} &  A_{22}   }\left[\begin{array}{c} u\\   -v \end{array}\right] 
 =(3-t)\left[\begin{array}{c} u\\   -v \end{array}\right]  $$
where $A_{11}$ is $(2 \times 2)$, $A_{22}$ is $(3 \times3)$, and $z=[u \ -v]^t$ where $u>0$ and $v\geq 0$ are vectors in $\mathbb{R}^2$ and $\mathbb{R}^3$ respectively.\\
Thus $A_{11}u-A_{12}v=(3-t)u$ implies $(A_{11}-(3-t)I_2)u \geq 0$ since the vector $A_{12}v$ is non-positive. But this implies, by a corollary of the Courant-Fischer theorem, (see Theorem 8.3.2 of \cite{Horn}), that $\rho(A_{11})\geq (3-t)$, where $\rho$ denotes the Perron root of $A_{11}$. Thus 
$A_{11}$ has an eigenvalue $\nu$ with
\begin{eqnarray}\label{CF}
 \nu \geq (3-t). 
\end{eqnarray}
Now $A_{11}$ is a $(2 \times 2)$ nonnegative matrix with trace zero since $tr(A)=0$ and so its eigenvalues are say, $\nu$ and $-\nu$. But $A$ is symmetric which means the interlacing property holds i.e. any principal submatrix of $A$ cannot have an eigenvalue larger than the largest eigenvalue of $A$ or an eigenvalue smaller than the smallest eigenvalue of $A$. \\
Hence the eigenvalues of $A_{11}$ must satisfy $-2\leq \nu,-\nu \leq 3+t$.\\
By \eqref{CF}, this means $-2 \leq -\nu \leq t-3$ which gives $t\geq 1$. \\
This argument is due to McDonald and Neumann \cite{McD/Neu} and we will make use of it again later when the condition \emph{symmetric} is replaced by \emph{diagonalizable}.\\
More generally, for a list of five real numbers $\lambda_1, \lambda_2, \lambda_3, \lambda_4 , \lambda_5$ satisfying \\
$\lambda_1\geq \lambda_2\geq \lambda_3 \geq \lambda_4 \geq \lambda_5$, McDonald and Neumann \cite{McD/Neu} showed that
$\lambda_2+\lambda_5\leq$ trace$(A)$ where $A$ is a symmetric realizing matrix.\\ Hence for $\widehat{\sigma}_t$ we must have that $(3-t)+(-2)\leq 0$ which again implies $t\geq 1$.\\
For $t=1$, $\widehat{\sigma}_t=(3+t,3-t,-2,-2,-2)$ is symmetrically realizable by
$$A=\npmatrix{ 0 & 2 & 0 & 0 & 0  \cr 2 & 0 & 0 & 0 & 0 \cr
0 & 0 & 0 & 2 & 2 \cr 0 & 0 & 2 & 0 & 2 \cr 0 &
0 & 2 & 2 & 0 \cr}$$
and so $t=1$ is the best possible result in the symmetric case - this was first proven by Loewy and Hartwig (unpublished).


\subsubsection{$\widehat{\sigma}_t=\{3+t,3-t,-2,-2,-2\}$ in the general case}
However in the general (non-symmetric) case Laffey and Meehan show that $\widehat{\sigma_t}$ must satisfy the necessary condition 
$$4s_4\geq s_2^2.$$ This requires that\\
$$t\geq t_0:=\sqrt{16\sqrt{6}-39}=0.43799\cdots$$ and they show that the matrix
$$A=\npmatrix{ 0 & 1 & 0 & 0 & 0  \cr \frac{15+t^2}{2} & 0 & 1 & 0 & 0 \cr
0 & 0 & 0 & 1 & 0 \cr 0 & 0 & 0 & 0 & 1 \cr 3t^4+58t^2+3 & \frac{t^4+78t^2-15}{4} &
10+6t^2 & \frac{15+t^2}{2} & 0  \cr}$$
is nonnegative for $t\geq t_0$ and has spectrum $(3+t,3-t,-2,-2,-2)$.\\
This result shows that the belief that $t=1$ was the minimum $t>0$ required for realization of $\widehat{\sigma}_t$ in the general NIEP is false. \\
Also note, that by the result of Guo \cite{Guo}, the list $(3+\widehat{t},3,-2,-2,-2)$ is realizable \\
for all $\widehat{t}\geq 2t_0=0.87598\cdots$.\\

\section{D-RNIEP $\neq$ SNIEP}
Next we examine the subtle difference between the SNIEP and the Diagonalizable RNIEP or D-RNIEP, where the D-RNIEP is the problem of finding a nonnegative \emph{diagonalizable} matrix realizing a given real spectrum $\sigma$.\\
Again we consider the list $\widehat{\sigma}_t=(3+t,3-t,-2,-2,-2)$.
We note that for 
$$t\geq t_0=\frac{1}{10}\sqrt{120\sqrt{1066}-3899}=0.4354153419\cdots$$
the perturbed list $(3+t,3-t,-1.9,-2,-2.1)$ is the spectrum of the matrix 
$$A=\npmatrix{ 0 & 1 & 0 & 0 & 0  \cr \frac{1501+100t^2}{200} & 0 & 1 & 0 & 0 \cr
0 & 0 & 0 & 1 & 0 \cr 0 & 0 & 0 & 0 & 1 \cr \frac{15000t^4+289950t^2+14649}{5000} &
\frac{10000t^4+779800t^2-148199}{40000} & \frac{150t^2+249}{25} & \frac{1501+100t^2}{200} & 0 \cr}.$$
The $(5,2)$ entry of the matrix A is nonnegative if $10000t^4+779800t^2\geq 148199$ and this holds for $t\geq t_0$.
Also note that this matrix is \emph{diagonalizable} as it has five distinct eigenvalues.\\
However if the list $(3+t,3-t,-1.9,-2,-2.1)$ is to be symmetrically realizable by a matrix $A$ we can use the Courant-Fischer argument again to say that the Perron root $\nu_1$ of some principal $2 \times 2$ submatrix $A_{11}$ say, must be at least $3-t$.\\ Let $\nu_1=3-t+\epsilon$ where $\epsilon \geq0$.\\ The eigenvalues $\nu_1,\nu_2$ of $A_{11}$ must satisfy $\nu_1+\nu_2=0$ (since tr$(A_{11}=0$).\\ Now consider the matrix $A+(2.1)I_5$.\\
This matrix has eigenvalues $5.1+t,5.1-t,0.2,0.1$ and $0$ and so is positive-semidefinite.
Now tr$(A_{11}+(2.1)I_2)=4.2$ since tr$(A_{11})=0$, and from this we get that,
$$4.2= 5.1-t+\epsilon+\mu_2$$
where $\mu_2=\nu_2+2.1$.
Since every principal submatrix of a positive-semidefinite matrix inherits positive-semidefiniteness we must have that $\mu_2\geq0$ and so we get 
$$\mu_2=-0.9+t-\epsilon \geq 0 \implies t\geq 0.9.$$
Hence the two problems of D-RNIEP and SNIEP are different at least in this case.\\

We also note that for small $\epsilon>0$ the spectrum $(3+t,3-t,-2+\epsilon,-2,-2-\epsilon)$ is diagonalizably realizable (five distinct eigenvalues) for $t$ close to $\sqrt{16\sqrt{6}-39}=0.43799\cdots$.\\
Spector showed that for this $t$ the same list is symmetrically realizable also.
However this is not a continuous property in $\epsilon$ since we have the following
\begin{Pro}\label{tprop}
Suppose $\widehat{\sigma}_t=(3+t,3-t,-2,-2,-2)$ is realizable by a diagonalizable matrix $A$, then $t\geq 1$.
\end{Pro}
To prove this result we will make use of the following fact.
\begin{Lem}\label{Schur}
Let $$B=\npmatrix{ B_{11}  & B_{12}\cr   
 B_{21} &  B_{22} }$$ be an $n\times n$ matrix.\\
If rank$(B)= $rank$(B_{11})=k$ and $B_{11}^{-1}$ exists, then $B_{22}=B_{21}B_{11}^{-1}B_{12}.$
\begin{proof}
Note that
$$\npmatrix{ I  & 0 \cr   
 -B_{21}B_{11}^{-1} &  I }\npmatrix{ B_{11}  & B_{12}\cr   
 B_{21} &  B_{22} }=\npmatrix{ B_{11}  & B_{12}\cr   
 0 &  B_{22} - B_{21}B_{11}^{-1}B_{12}},$$
 where $I$ and $0$ are the identity and zero matrices of appropriate dimensions respectively.
Since rank$(B)$=rank$(B_{11})=k$ the rowspace of $B$ is spanned by the rows of $B_{11}$ extended in to $B_{12}$. Hence for any row $r_j,$ $j=k+1,k+2,\ldots,n,$ \\ 
in $B_{22} - B_{21}B_{11}^{-1}B_{12}$ to be written as a linear combination of rows in $B_{11}$ extended we must have that \\ $r_j=\alpha_1r_1+\alpha_2r_2+\cdots+\alpha_kr_k+\alpha_{k+1}r_{k+1}+\cdots+\alpha_nr_n$.
But $r_1,r_2,\ldots,r_k$ are already linearly independent, so $\alpha_1=\alpha_2=\cdots = \alpha_n=0$. Hence $r_j=0$ for $j=k+1,k+2,\ldots, n$.
\end{proof}
\end{Lem}

\noindent $\bf{Proof \ of \ Proposition \ \ref{tprop}}$
\begin{proof}
Suppose $t<1$.\\
Note for $\widehat{\sigma}_t=(3+t,3-t,-2,-2,-2)$ to be realizable by a diagonalizable matrix $A$, then $A$ must be irreducible.\\
Let $w^T>0$ be the left eigenvector associated with the Perron eigenvalue $3+t$, so $w^TA=(3+t)w^T$, where $T$ denotes the transpose.\\
Let $v$ be the right eigenvector associated with the eigenvalue $3-t$, so $Av=(3-t)v$.\\
Then $(3+t)w^Tv=w^TAv=w^T(3-t)v \Rightarrow w^Tv=0 \Rightarrow v$ has some negative components.
By the earlier argument we know that $v$ cannot have just one negative (or positive) component.\\
Hence we can write $Az=(3-t)z$ as
$$\npmatrix{ A_{11}  & A_{12}\cr   
 A_{21} &  A_{22}   }\left[\begin{array}{c} z_1\\   -z_2 \end{array}\right] 
 =(3-t)\left[\begin{array}{c} z_1\\   -z_2 \end{array}\right]  $$
where $A_{11}$ is ($2\times 2$), $A_{22}$ is ($3 \times3)$,   and $z_1,z_2 \geq 0$ are column vectors in $\mathbb{R}^2$ and $\mathbb{R}^3$ respectively.
Then (again using the McDonald-Neumann argument) the Perron eigenvalue of $A_{11}$ is at least $3-t$ and hence $A_{11}$ has spectrum $\pm(3-t+\epsilon)$ where $\epsilon \geq0$, since tr$(A_{11})=0$.
Since $A$ is diagonalizable the minimum polynomial of $A$ is 
$$m_A(x)=(x-(3+t))(x-(3-t))(x+2)$$ so that
\begin{eqnarray*}\label{trace}
&&(A-(3+t)I)(A-(3-t)I)(A+2I)=0   \\
&\Rightarrow& (A^2-6A+(9-t^2)I)(A+2I)=0  \\
&\Rightarrow& A^3-4A^2-(3+t^2)A+(18-2t^2)I=0  
\end{eqnarray*}
\begin{eqnarray}
&\Rightarrow& A^3+(18-2t^2)I=4A^2+(3+t^2)A. 
\end{eqnarray}
Now $A$ diagonalizable also implies there exists a nonsingular matrix $T$ with
$$T^{-1}(A+2I_5)T=(5+t)\oplus(5-t)\oplus(0_3).$$
Now rank$(A+2I)=2$ since rank$(A_{11}+2I_2)=2.$\\
Note that $A_{11}$ has trace zero and so $A_{11}+2I$ has the form
$$\npmatrix{ 2  & a_{12}\cr   
 a_{21} &  2 }. $$
So det$(A_{11}+2I)=4-a_{12}a_{21}$

Since rank$(A+2I)=2$ we have that rank$(A_{22}+2I_3)\leq 2$ since the rank of a leading principal submatrix cannot exceed the rank of the original matrix. Thus $A_{22}+2I_3$ has zero as an eigenvalue (as it does not have full rank) and hence $-2$ is an eigenvalue of $A_{22}$. \\
Applying the McDonald-Neumann argument again we have that $A_{21}z_1-A_{22}z_2=-(3-t)z_2$ so that $(A_{22}-(3-t))z_2=A_{21}z_1\geq 0$ and so $A_{22}$ has an eigenvalue $3-t+\epsilon'$ where $\epsilon' \geq0$.\\
Since tr$(A_{22})=0$, $A_{22}$ has eigenvalues $3-t+\epsilon',-1+t-\epsilon',-2$.\\
Now consider the tr$(A^2)$ where
\begin{eqnarray*}A^2=
     \left(\begin{array}{cc}
      A_{11}^2+A_{12}A_{21} & A_{11}A_{12}+A_{12}A_{22} \\
      A_{21}A_{11}+A_{22}A_{21} & A_{21}A_{12}+A_{22}^2 \\
    \end{array}\right) &  & \\
    &  & 
  \end{eqnarray*}
  
and the the tr$(A^3)$ where $A^3=$
     \begin{eqnarray*}=\left(\begin{array}{cc}
      A_{11} & A_{12} \\
      A_{21} & A_{22} \\
    \end{array}\right) 
     \left(\begin{array}{cc}
      A_{11}^2+A_{12}A_{21} & A_{11}A_{12}+A_{12}A_{22} \\
      A_{21}A_{11}+A_{22}A_{21} & A_{21}A_{12}+A_{22}^2 \\
    \end{array}\right) 
  \end{eqnarray*}
\begin{eqnarray*}=
     \left(\begin{array}{cc}
      A_{11}^3+A_{11}A_{12}A_{21}+A_{12}\left(A_{21}A_{11}+A_{22}A_{21}\right) & * \\
      ** & A_{21}\left(A_{11}A_{12}+A_{12}A_{22}\right)+A_{22}^3+A_{22}A_{21}A_{12} \\
    \end{array}\right) &  & \\
    &  & 
  \end{eqnarray*}
where * and ** do not contribute to the trace of $A^3$.\\  
Let $$\alpha= \mbox{tr}(A_{12}A_{21})=\mbox{tr}(A_{21}A_{12})$$
$$\beta = \mbox{tr}(A_{11}A_{12}A_{21}) \mbox{ \ and}$$
$$\gamma = \mbox{tr}(A_{12}A_{22}A_{21}).$$
Using \eqref{trace} we see that the contribution to tr$I$, tr$(A)$, tr($A^2)$ and tr($A^3)$ from positions $(1,1)$ and $(2,2)$ yields that
$$4\left(2(3-t+\epsilon)^2+\alpha\right)=2\beta+\gamma+2(18-2t^2).$$

Applying Lemma 2 to $A+2I$=
$$\npmatrix{ A_{11}+2I_2  & A_{12}\cr   
 A_{21} &  A_{22}+2I_3   }$$
 we get that
 $A_{22}+2I_3=A_{21}(A_{11}+2I_2)^{-1}A_{12}$.\\
 Now $(A_{11}+2I_2)^{-1}$
 $$=\frac{1}{\mbox{det}(A_{11}+2I_2)}\npmatrix{ 2  & -a_{12}\cr   
 -a_{21} &  2 }$$
 $$=\frac{-1}{(5-t+\epsilon)(1-t+\epsilon)}\npmatrix{ 2  & -a_{12}\cr   
 -a_{21} &  2 }$$
 $$=\frac{1}{(5-t+\epsilon)(1-t+\epsilon))}\npmatrix{ -2  & a_{12}\cr   
 a_{21} &  -2 }$$
 $$=\frac{(A_{11}-2I_2)}{(5-t+\epsilon)(1-t+\epsilon))}.$$
 So $A_{22}+2I_3=\frac{A_{21}(A_{11}-2I_2)A_{12}}{(5-t+\epsilon)(1-t+\epsilon))}$.
 Note that tr$A_{22}=0$ so upon comparing traces in this equation we get that\\
 $$6=\frac{\mbox{tr}(A_{21}A_{11}A_{12})-2\mbox{tr}(A_{21}A_{12})}{(5-t+\epsilon)(1-t+\epsilon))}$$
 $$=\frac{\beta-2\alpha}{(5-t+\epsilon)(1-t+\epsilon))}.$$
Now consider the contribution from the top two positions on the main diagonal to the traces of both sides of equation \eqref{trace}. We have that  
 $$2\beta+\gamma+36-4t^2=8(3-t+\epsilon)^2+4\alpha$$
This implies
$$12(5-t+\epsilon)(1-t+\epsilon)+4\alpha+\gamma+36-4t^2=72-48t+8t^2+16\epsilon(3-t)+8\epsilon^2$$
$\Rightarrow$
$$24+24\epsilon+4\epsilon^2+4\epsilon+\gamma-(24t+8\epsilon t)=0.$$
But note that $24-24t>0$ since $t<1$, and $24\epsilon-8\epsilon t=16\epsilon +8 \epsilon-8\epsilon t >0$, implies that the left hand side of this equation is positive and the right hand side is zero which contradicts our hypothesis that $t<1$. Hence $t\geq1.$
\end{proof}
\vspace{1cm}
Laffey and Smigoc \cite{LafSmi2} proved that $(3+t,3-t,-2,-2,-2,0)$ is symmetrically realizable by a $6 \times 6$ matrix for $t\geq \frac{1}{3}$. The fact that the value $t=\frac{1}{3}$ is necessary in this case and in the case $\widehat{\sigma}_t\bigcup(0,0)$ has been checked by computer search by Lixing Han (unpublished) at the University of Michigan-Flint and the same realizing matrix was found.  
We conjecture that $t=\frac{1}{3}$ is the best bound for any number of zeros added to the spectrum $\widehat{\sigma}_t$. Also Lixing Han's extensive computer searches in the cases $n = 6$ and $n = 7$ failed to find a counterexample to this bound.

\section{A note on the dependence of realizable spectra on the Jordan Form structure}
The matrix
$$A=\frac{1}{4}\npmatrix{ 0 & 8 & 1 & 0 & 0  \cr 8 & 0 & 1 & 0 & 0 \cr
\frac{75}{2} & \frac{75}{2} & 0 & 1 & 0 \cr 0 & 0 & 0 & 0 & 1 \cr 829 &
829 & 256 & 110 & 0 \cr}$$
has spectrum $\sigma_{\frac{3}{4}}=(3+\frac{3}{4},3-\frac{3}{4},-2,-2,-2)$ and minimal polynomial
$$\left(x-(3+\frac{3}{4}))(x-(3-\frac{3}{4}))(x+2)^2\right)$$
and so it distinguishes the diagonalizably realizable lists from those with Canonical Jordan Form structure $(3+t)\oplus(3-t)\oplus(-2)\oplus\npmatrix{ -2 & 1 \cr 0 & -2 \cr}$ with $t=\frac{3}{4}$.\\
Hence the D-RNIEP is different to the general NIEP.
Note that this matrix was built up from a $4\times 4$ nonnegative matrix with spectrum $(3+\frac{3}{4},3-\frac{3}{4},-2,-2)$ having the entry 2 on the diagonal using the \v{S}migoc methods deployed in \cite{Smi}. The question of whether every realizing matrix can be chosen to be nonderogatory remains open and will require further ideas related to this paper.

\section{Conclusion}
In this paper we proved that the SNIEP $\neq$ D-RNIEP and that the D-RNIEP can be distinguished from the general NIEP by examining the Jordan Normal Form. We also showed that the minimum $t>0$ for which $(3+t,3-t,-2,-2,-2)$ is realizable by a diagonalizable matrix is $t=1$.

\end{document}